\documentclass[12pt]{article}

\usepackage{setspace}

\usepackage{amsthm, amssymb, amstext}
\usepackage[fleqn]{amsmath}
\usepackage{latexsym}
\usepackage[dvips]{graphicx}
\usepackage{comment}
\usepackage{hyperref}
\usepackage[capitalize]{cleveref}
\usepackage{mathtools}
\usepackage{enumerate} 
\usepackage{paralist}
\usepackage{enumitem}

\usepackage{todonotes}
\usepackage{comment}

\crefname{claim}{Claim}{Claims}

\newtheorem{theorem}{Theorem}[section]
\newtheorem{lemma}[theorem]{Lemma}

\newtheorem{claim}{Claim}

\theoremstyle{definition}

\usepackage{amsthm}

\definecolor{myRed}{rgb}{0.68, 0.05, 0.0}
\colorlet{myBlue}{blue!70!black}
\colorlet{myViolet}{myBlue!55!myRed}
\definecolor{darkraspberry}{rgb}{0.53, 0.15, 0.34}
\definecolor{olive}{rgb}{0.42, 0.56, 0.14}

\hypersetup{
	colorlinks=true,
        linkcolor=myBlue, 
        citecolor=myBlue,
	bookmarksopen=true,
	bookmarksnumbered,
	bookmarksopenlevel=2,
	bookmarksdepth=3
}

\usepackage{authblk}

\title{Graphs without a $3$-connected subgraph are $4$-colourable}

\author[1]{\'Edouard Bonnet}
\author[1]{Carl Feghali}
\author[2]{Tung Nguyen}
\author[3]{Alex~Scott}
\author[2]{Paul~Seymour}
\author[1]{St\'ephan Thomass\'e}
\author[1]{Nicolas Trotignon}

\affil[1]{Univ Lyon, EnsL, CNRS, LIP, F-69342, Lyon Cedex 07, France}
\affil[2]{Princeton University, Princeton, NJ 08544, USA}
\affil[3]{Mathematical Institute, University of Oxford, Oxford OX2 6GG, UK}

\date{\today}

\begin{document}
\maketitle

\begin{abstract}
  In 1972, Mader showed that every graph without a $3$-connected
  subgraph is \mbox{$4$-degenerate} and thus \mbox{$5$-colourable}. We show
  that the number $5$ of colours can be replaced by $4$, which is best
  possible.
 \end{abstract}

\section{Introduction}

Throughout the paper all graphs are finite and simple, and we only use standard notions and notation.  
We recall that a~graph is \emph{$k$-connected} if it has at~least~$k+1$ vertices and no vertex cutset with at~most~$k-1$ vertices. 
In 1972, Mader~\cite{mader:mkc} proved the following theorem.  

\begin{theorem}
  \label{th:mader}
  For every integer $k\geq 1$, every graph with average degree at least
  $4k$ contains a $(k + 1)$-connected subgraph.
\end{theorem}

\begin{figure}
    \begin{center}
    \begin{tabular}{ccc}
    \includegraphics[height=3cm]{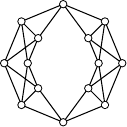}  & \rule{1cm}{0cm} 
    &     \includegraphics[height=3cm]{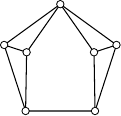}\\
    Minimum degree 4 & & Chromatic number 4      
    \end{tabular}
  \end{center}
  \caption{Graphs with no 3-connected
    subgraph.
\label{f:4}\label{f:wheel}}
  \end{figure}

Focusing on the case $k = 2$ of~\cref{th:mader}, we call a~graph \emph{fragile} if it has no \mbox{3-connected} subgraph. From~\cref{th:mader}, every non-null fragile graph has a vertex of degree at most~7.  
By restricting the proof of Mader to the case $k=2$, it is easy to show that all fragile graphs $G$ on at least four vertices satisfy $|E(G)|\leq 2.5|V(G)| -5$ (we supply the proof in Section~\ref{sec:c} for the sake of completeness).  So the average degree of $G$ is smaller than $5$.  Thus every fragile graph contains a vertex of degree at most~4, and this is best possible as shown by the graph in Figure~\ref{f:4}. Every fragile graph is therefore 5-colourable. 

Despite recent progress on related questions, there is no available proof that the number $5$ of colours can be improved. The objective of this paper is to prove the following theorem that implies that every fragile graph is 4-colourable. It was announced without proof in~\cite{nguyen:conn} (which also contains a thorough literature review) and was independently rediscovered by the first two and last two authors of this article.

\begin{theorem}\label{thm:main}
For all $m\geq 4$, every graph with chromatic number at least~$m+1$ has a 3-connected subgraph with chromatic number at least $m$.
\end{theorem}

Theorem~\ref{thm:main} is best possible as shown by the graph
in Figure~\ref{f:wheel}.  The proof of
Theorem~\ref{thm:main} is given in Section~\ref{sec:p}.  Several
remarks and open questions are presented in Section~\ref{sec:c}.

\section{Proof of Theorem \ref{thm:main}}
\label{sec:p}

For every integer $m\geq 4$, a graph $G$ is \emph{$m$-fragile} if all 3-connected subgraphs of $G$ are $(m-1)$-colourable. Observe that a fragile graph is $m$-fragile for all $m\geq 4$. \cref{thm:main} can be rephrased as: for all $m\geq 4$, every $m$-fragile graph is $m$-colourable. To prove \cref{thm:main}, we shall establish the following stronger statement. By a \emph{$k$-colouring of a graph $G$}, we mean a function $c$ that associates to each vertex of $G$ an integer in $\{1, \dots, k\}$ and such that for all edges $xy$ of $G$, $c(x)\neq c(y)$. 

\begin{theorem}\label{thm:main2}
For every integer $m\geq 4$, every $m$-fragile graph $G$ satisfies the following four
 conditions.
\begin{enumerate}
\item\label{c1} For all non-adjacent $x, y\in V(G)$, $G$ admits an
  $m$-colouring $c$ such that $c(x) = c(y)$.
\item\label{c2} For all distinct $x, y \in V(G)$, $G$ admits an $m$-colouring $c$
  such that $c(x) \neq c(y)$.
\item\label{c3} For all distinct $x, y, z \in V(G)$, $G$ admits an $m$-colouring $c$ such that
  $c(x) \notin \{c(y), c(z)\}$.
\item\label{c4} For all distinct $x, y, z \in V(G)$ that are not
  all pairwise adjacent, $G$ admits an $m$-colouring $c$ such
  that $|\{c(x), c(y), c(z)\}| = 2$.
\end{enumerate}
\end{theorem}

\begin{proof}
We proceed by induction on $|V(G)|$. If $|V(G)| \leq 3$, then $G$ obviously satisfies conditions~\ref{c1}--\ref{c4}. For the induction step, suppose $|V(G)|\geq 4$ and that the statement holds for every graph with fewer vertices than $G$.  

If $G$ is 3-connected, then $G$ satisfies conditions~\ref{c1}--\ref{c4} because by assumption $G$ is  $(m-1)$-colourable. So $G$ can be coloured with colours 1 to $m-1$, and colour $m$  is available to satisfy any of the conditions~\ref{c1}--\ref{c4} (for instance $x$ and $y$ can be recoloured with colour $m$ to satisfy~\ref{c1}). Hence we may assume from here on that $G$ is not 3-connected.

Since $G$ is not 3-connected, there exist two induced subgraphs $G_1, G_2$ of $G$
such that $V(G) = V(G_1) \cup V(G_2)$, $E(G) = E(G_1) \cup E(G_2)$,
$V(G_1) \setminus V(G_2) \neq \emptyset$,
$V(G_2) \setminus V(G_1) \neq \emptyset$, and
$S = V(G_1) \cap V(G_2)$ has size at~most~2.  Moreover, since $G$ is $m$-fragile, $G_1$ and $G_2$ are also
$m$-fragile and, as $|V(G_1)|, |V(G_2)| < |V(G)|$, we may apply the
induction hypothesis to both $G_1$ and~$G_2$.

If $S= \emptyset$, the induction step is obvious and we omit the
details. So we may set $S = \{u, v\}$ (possibly $u =v$). We have to prove that for each of the precolouring 
conditions~$C$ among \ref{c1}--\ref{c4} on any given set $X \subseteq V(G)$ (namely, $X=\{x, y\}$ for conditions \ref{c1} and \ref{c2} and $X=\{x, y, z\}$ for conditions \ref{c3} and
\ref{c4}) some appropriate 4-colouring exists.  Suppose first that $X\subseteq V(G_1)$. Then, by the induction hypothesis, $G_1$ admits a colouring $c_1$ that satisfies~$C$.  By applying~\ref{c1} or~\ref{c2} to the vertices $u$ and $v$ of $G_2$ (or trivially if $u=v$), and up to a relabeling of the colours, we can force a colouring $c_2$ of $G_2$ such that $c_2(u) = c_1(u)$ and $c_2(v) = c_1(v)$. Note that the case when $uv$ is an edge corresponds to the usual amalgamation of two colourings on a clique cutset. Hence, $c_1 \cup c_2$ is a colouring of $G$ that satisfies~$C$.  The proof is similar
when $X \subseteq V(G_2)$.  Hence, from here on, we may assume that
\[ \text{$X$ intersects both $V(G_1) \setminus V(G_2)$ and
$V(G_2) \setminus V(G_1)$.
}\tag{$\star$} \]

We now prove four claims, from which Theorem \ref{thm:main2} trivially follows.  Their proofs are easy  when $u=v$, so we omit this case and assume $u\neq v$.  Note that, unless specified otherwise, we shall
make no assumption on whether $u$ and $v$ are adjacent.

\begin{claim}
  \label{l:c1}
  The graph $G$ satisfies \ref{c1}.
\end{claim}

\begin{proof}
  By ($\star$), we may assume that $x\in V(G_1)\setminus V(G_2)$ and
  $y\in V(G_2)\setminus V(G_1)$.  We build three colourings $a_1$,
  $b_1$ and $c_1$ of $G_1$ and three colourings $a_2$, $b_2$ and $c_2$
  of $G_2$ that are represented in Figure~\ref{f:c1} for the reader's
  convenience.

\begin{figure}
    \begin{center}
    \begin{tabular}{cccc}
    & $G_1$ &\parbox{2cm}{}& $G_2$\\
    &
    \rule{0cm}{1.5cm}\begin{minipage}[c][2.3cm]{2cm}\includegraphics[height=2cm]{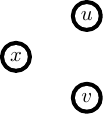}\end{minipage}&&
    \begin{minipage}[c][2.3cm]{2cm}\includegraphics[height=2cm]{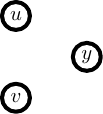}\end{minipage}\\\hline
    \begin{minipage}[c][2cm]{.5cm}$a$\end{minipage}&
    \begin{minipage}[c][2.3cm]{2cm}\includegraphics[height=2cm]{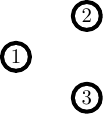}\end{minipage}&&
    \begin{minipage}[c][2.3cm]{2cm}\includegraphics[height=2cm]{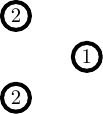}\end{minipage}\\\hline
    \begin{minipage}[c][2cm]{.5cm}$b$\end{minipage}&
    \begin{minipage}[c][2.3cm]{2cm}\includegraphics[height=2cm]{graphCut-8.pdf}\end{minipage}&&
    \begin{minipage}[c][2.3cm]{2cm}\includegraphics[height=2cm]{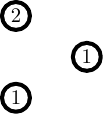}\end{minipage}\\\hline
    \begin{minipage}[c][2cm]{.5cm}$c$\end{minipage}&
    \begin{minipage}[c][2.3cm]{2cm}\includegraphics[height=2cm]{graphCut-8.pdf}\end{minipage}&&
    \begin{minipage}[c][2.3cm]{2cm}\includegraphics[height=2cm]{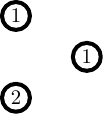}\end{minipage}\\
\end{tabular}
  \end{center}
  \caption{Colourings obtained in the proof of~\cref{l:c1}.
\label{f:c1}}
  \end{figure}

  By~\ref{c3} applied to $x, u, v$ (in this order) in $G_1$, we obtain a colouring
  $a_1$ of $G_1$ such that $a_1(x) \notin \{a_2(u), a_2(v)\}$.
  Similarly, we obtain a colouring $a_2$ of $G_2$ such that
  $a_2(y) \notin \{a_2(u), a_2(v)\}$.  Up to a relabeling, we may
  assume that $a_1(x) = a_2(y) = 1$, $a_1(u) = a_2(u) = 2$ and
  $a_1(v), a_2(v) \in \{2, 3\}$.  If $a_1(v) = a_2(v)$, then
  $a_1 \cup a_2$ is a colouring of $G$ that satisfies~\ref{c1}.  Hence,
  up to symmetry, we may assume that $a_1(v) = 3$ and $a_2(v) = 2$.

  By~\ref{c3} applied  to $u, x, v$ in $G_1$, we obtain a colouring
  $b_1$ of $G_1$ such that $b_1(u) \notin \{b_1(x), b_1(v)\}$.
  Similarly, we obtain a colouring $b_2$ of $G_2$ such that
  $b_2(u) \notin \{b_2(y), b_2(v)\}$.  Up to a relabeling, we may
  assume that $b_1(x) = b_2(y) = 1$, $b_1(u) = b_2(u) = 2$ and
  $b_1(v), b_2(v) \in \{1, 3\}$.  If $b_1(v) = b_2(v)$, then
  $b_1 \cup b_2$ is a colouring of $G$ that satisfies~\ref{c1}.  Hence,
  we may assume that $b_1(v) \neq b_2(v)$.  If $b_2(v) = 3$, then
  $a_1 \cup b_2$ is a colouring of $G$ that satisfies~\ref{c1}.  Hence,
  we may assume that $b_1(v) = 3$ and $b_2(v) = 1$.

  By~\ref{c3} applied to $v, x, u$ in $G_1$, we obtain a colouring
  $c_1$ of $G_1$ such that $c_1(v) \notin \{c_1(x), c_1(u)\}$.
  Similarly, we obtain a colouring $c_2$ of $G_2$ such that
  $c_2(v) \notin \{c_2(y), c_2(u)\}$.  Up to a relabeling, we may
  assume that $c_1(x) = 1$ and either $c_1(u) = 1 $ and $ c_1(v) = 2$ 
  or $c_1(u) = 2$ and $ c_1(v) = 3$.  Up to a relabeling, we may also
  assume that $c_2(y) = 1$ and either $c_2(u) = 1 $ and $ c_2(v) = 2$
  or $c_2(u) = 2$ and $ c_2(v) = 3$.  If $c_1(u) = c_2(u)$, then
  $c_1 \cup c_2$ is a colouring of $G$ that satisfies~\ref{c1}.  Hence,
  we may assume that $c_1(u) \neq c_2(u)$.  If $c_2(u) = 2$ (and so
  $c_2(v) = 3$), then $a_1 \cup c_2$ is a colouring of $G$ that
  satisfies~\ref{c1}.  Hence, we may assume that $c_1(u) = 2$,
  $c_1(v) = 3$, $c_2(u) = 1$ and $c_2(v) = 2$.

  By~\ref{c4} applied to $x, u, v$ in $G_1$, we obtain a colouring
  $d_1$ of $G_1$ such that $|\{d_1(x), d_1(u), d_1(v)\}| = 2$ (note that $x$, $u$ and $v$ are not pairwise adjacent because $a_2(u)=a_2(v)$ implies $uv\notin E(G)$). Up to a
  relabeling, we may assume that $d_1(x) = 1$ and
  $\{d_1(x), d_1(u), d_1(v)\} = \{1, 2\}$.  If $d_1(u) = 1$ and $d_1(v) = 2$,
  then $d_1 \cup c_2$ satisfies~\ref{c1}. And  if $d_1(u) = 2$ and
  $d_1(v) = 1$, then $d_1 \cup b_2$ satisfies~\ref{c1}.  Finally, if
  $d_1(u) = 2$ and $d_1(v) = 2$, then $d_1 \cup a_2$
  satisfies~\ref{c1}. The claim is proved.
\end{proof}

\begin{claim}
  \label{l:c3}
  The graph $G$ satisfies \ref{c3}.
\end{claim}

\begin{proof}
  If $x\in \{u, v\}$ (say $x=u$ up to symmetry), then by ($\star$) we
  may assume that $y\in V(G_1) \setminus V(G_2)$ and
  $z \in V(G_2) \setminus V(G_1)$.   By~\ref{c3}
  applied separately to $x$, $v$ and $y$ in $G_1$ and to $x$, $v$ and $z$ in
  $G_2$, we obtain up to a relabeling a colouring $a_1$ of $G_1$ and a
  colouring $a_2$ of $G_2$ such that $a_1(x)=a_2(x) = 1$,
  $a_1(v)=a_2(v) = 2$, $a_1(y)\neq 1$ and $a_2(z)\neq 1$.  Hence,
  $a_1\cup a_2$ is a colouring of $G$ that satisfies $\ref{c3}$.  We may therefore assume that
  $x\notin \{u, v\}$, and so up to symmetry that
  $x\in V(G_1) \setminus V(G_2)$.
   
  Hence, by ($\star$) and up to symmetry, we may restrict our
  attention to the following two cases.

  \medskip
   
  {\noindent\bf Case 1:} $x\in V(G_1) \setminus V(G_2)$ and
  $y, z \in V(G_2)$.

    If $uv\in E(G)$, then by~\ref{c3} applied to $x$, $u$ and $v$ and up to a relabeling, there exists a colouring $a_1$ of $G_1$ such that  $a_1(x) = 1$, $a_1(u)=2$, and $a_1(v)=3$. 
    We claim that there exists a colouring $a_2$ of $G_2$ that requires at most $m-1$ colours for $u$, $v$, $y$, $z$.  
    If $m\geq 5$, this is trivial, so suppose $m=4$. Then the graph induced by $u$, $v$, $y$ and $z$ is not a complete graph on four vertices, because such a graph is 3-connected with chromatic number~4 and would imply that $G$ is not $m$-fragile. 
    Hence, either $|\{u, v, y, z\}|\leq 3$ or there are non-adjacent vertices among $u$, $v$, $y$ and $z$.  
    In either case, there exists a colouring $a_2$ of $G_2$ that requires at most $m-1=3$ colours for $u$, $v$, $y$, $z$ (trivially if $|\{u, v, y, z\}|\leq 3$ or by applying~\ref{c1} to a non-edge otherwise).  
    This proves our claim. 
    Up to a relabeling, we may assume that $a_2(u)=2$, $a_2(v)=3$ and $\{a_2(y), a_2(z)\} \subseteq \{2, \dots, m\}$.  Hence, $a_1\cup a_2$ is a colouring of $G$ satisfying~\ref{c3}. 
    We may therefore assume from here on that $uv\notin E(G)$. 

  Suppose that there exists a colouring $a_1$ of $G_1$ such that
  $a_1(x) \neq a_1(u) = a_1(v)$. So, up to a relabeling, we may assume $a_1(x)=1$ and $a_1(u)=a_1(v)=2$.  Then
  by~\ref{c1} applied to $u$ and $v$ in $G_2$, there exists a colouring
  $a_2$ of $G_2$ such that $a_2(u)=a_2(v)$.  Hence,
  $|\{a_2(u), a_2(v), a_2(y), a_2(z)\}|\leq 3$.  So, up to a
  relabeling, we may assume that $a_2(u) = a_2(v) = 2$  and $\{a_2(y), a_2(z)\} \subseteq \{2, 3, 4\}$.   So $a_1 \cup a_2$ is a colouring of $G$ that
  satisfies~\ref{c3}.  We may therefore assume that no colouring
  as $a_1$ exists.

  Hence, when applying~\ref{c3} to $x$, $u$ and $v$, up to a
  relabeling, we obtain a colouring $b_1$ of $G_1$ such that
  $b_1(x) =1$, $b_1(u)=2$ and $b_1(v)=3$. And when applying~\ref{c4}
  to $x$, $u$ and~$v$ (which is allowed since $uv\notin E(G)$), up to a relabeling and to the symmetry between
  $u$~and~$v$, we obtain a colouring $c_1$ of $G_1$ such that
  $c_1 (x)=1$, $c_1 (u)=1$ and $c_1 (v)=2$.  

  By~\ref{c2} applied to $u$ and $v$, there exists a colouring $d_2$ of
  $G_2$ such that $d_2(u)\neq d_2(v)$. If
  $|\{d_2(u), d_2(v), d_2(y), d_2(z)\}| \leq 3$, then up to a
  relabeling, we may assume that $d_2(u)=2$, $d_2(v)=3$ and
  $\{d_2(y), d_2(z)\}\subseteq \{2, 3, 4\}$, So
  $b_1 \cup d_2$ is a colouring of $G$ that satisfies~\ref{c3}.  And if
  $|\{d_2(u), d_2(v), d_2(y), d_2(z)\}| = 4$, then we may assume up to a
  relabeling that $d_2(u)=1$, $d_2(v)=2$, $d_2(y)=3$ and $d_2(z)=4$,
  so $c_1\cup d_2$ is a colouring that satisfies~\ref{c3}.

  \medskip
   
  {\noindent\bf Case 2:} $x, y\in V(G_1) \setminus V(G_2)$ and
  $z \in V(G_2) \setminus V(G_1)$.

  By~\ref{c3} applied to $x$, $y$ and $u$, up to a relabeling, we
  obtain a colouring $a_1$ of $G_1$ such that $a_1(x) =1$, $a_1(y)=2$
  and $a_1(u)\in \{2, 3\}$.  If $a_1(v)\neq 1$, then colour~1 is not
  used on $u$ or $v$ under $a_1$.  By~\ref{c1} or~\ref{c2} applied to $u$ and
  $v$, we obtain up to a relabeling a colouring $a_2$ of $G_2$ such
  that $a_2(u)=a_1(u)$ and $a_2(v)=a_1(v)$.  Thus,  colour 1 
  is not used on $u$ or $v$ under $a_2$ either and so, up to a relabeling, we may assume
  that $a_2(z) \neq 1$.  Hence $a_1 \cup a_2$ is a colouring of $G$ that
  satisfies~\ref{c3}.  We may therefore assume that $a_1(v)=1$.

  By~\ref{c3} applied to $v$, $u$ and $z$, up to a relabeling, we
  obtain a colouring $b_2$ of $G_2$ such that $b_2(v) =1$,
  $b_2(u)= a_1(u)$ and $b_2(z)\neq 1$.  Hence $a_1 \cup\,b_2$ is a
  colouring of $G$ that satisfies~\ref{c3}.
 \end{proof}

\begin{claim}
  \label{l:c2}
  The graph $G$ satisfies \ref{c2}.
\end{claim}

\begin{proof}
  By Claim~\ref{l:c3}, we may apply~\ref{c3} to $x$, $y$ and any vertex of $G$. We obtain a colouring of $G$ that satisfies~\ref{c2}.
\end{proof}

\begin{claim}
  \label{l:c4}
  The graph $G$ satisfies \ref{c4}.
\end{claim}

\begin{proof}
  By ($\star$), we may assume that $x\in V(G_1)\setminus V(G_2)$ and
  $y \in V(G_2)\setminus \{u\}$ and $z\in V(G_2) \setminus V(G_1)$.

    Suppose that $uv\in E(G)$. Then by~\ref{c3} applied to $x$, $u$ and $v$ and up to a relabeling, there exists a colouring $a_1$ of $G_1$ such that $a_1(x)=1$, $a_1(u)=2$ and $a_1(v)=3$. 
      By~\ref{c3} applied to $u$, $y$ and $z$ (that are distinct since $y\neq u$ and $z\in V(G_2) \setminus V(G_1)$) and up to a relabeling, we obtain a colouring $a_2$ of $G_2$ such that  $a_1(u)=2$, $a_1(v)=3$ and $\{a_2(y), a_2(z)\}$ is either $\{3, 1\}$, $\{3\}$ or $\{4\}$.  In either case, $a_1\cup a_2$
 is a colouring of $G$ satisfying~\ref{c4}.   
    We may therefore assume from here on that $uv\notin E(G)$. 
  
  Suppose that there exists a colouring $a_1$ of $G_1$ such that
  $a_1(x) \neq a_1(u) = a_1(v)$. Then up to a relabeling we may assume that $a_1(x)=1$ and $a_1(u)=a_1(v) = 2$. 
  By~\ref{c1} applied to $u$ and $v$ in $G_2$, there exists up to a relabeling a colouring
  $a_2$ of $G_2$ such that $a_2(u)=a_2(v) =2$.  If
  $a_2(y)=a_2(z)$, then up to relabeling, we may assume that
  $a_2(y)=a_2(z)\neq 1$, so~\ref{c4} is satisfied by $a_1 \cup a_2$. And if
  $a_2(y)\neq a_2(z)$, then up to a relabeling, we may assume
  $a_2(y) = 1$ or $a_2(z) = 1$, and~\ref{c4} is again
  satisfied by $a_1\cup a_2$.  We may therefore assume that no colouring as
  $a_1$ exists.

  Hence, when applying~\ref{c3} to $x$, $u$ and $v$, up to a
  relabeling, we obtain a colouring $b_1$ of $G_1$ such that
  $b_1(x) =1$, $b_1(u)=2$ and $b_1(v)=3$. And when applying~\ref{c4}
  to $x$, $u$ and $v$ (which is allowed since $uv\notin E(G)$), up to a relabeling and to the symmetry between
  $u$ and $v$, we obtain a colouring $c_1$ of $G_1$ such that
  $c_1 (x)=1$, $c_1 (u)=1$ and $c_1 (v)=2$.  

  On the other hand, by~\ref{c2} applied to $u$ and $v$, there exists a colouring $d_2$ of
  $G_2$ such that $d_2(u)\neq d_2(v)$.
If $d_2(y)= d_2(z)$, then up to a relabeling, we may assume that
  $d_2(u)=2$, $d_2(v)=3$ and $d_2(y) \neq 1$. Thus, $b_1\cup d_2$ is a
  colouring that satisfies~\ref{c4}.  Hence, from here on, we may assume that
  $d_2(y) \neq d_2(z)$.
  
  If $|\{d_2(u), d_2(v), d_2(y), d_2(z)\}| \geq 3$, then we may assume
  up to a relabeling that $d_2(u)=2$, $d_2(v)=3$ and
  $1\in \{d_2(y), d_2(z)\}$, so $b_1\cup d_2$ is a colouring that
  satisfies~\ref{c4}.  If $|\{d_2(u), d_2(v), d_2(y), d_2(z)\}| = 2$,
  then up to a relabeling, we may assume that $d_2(u)=1$, $d_2(v)=2$,
  so that $\{d_2(y), d_2(z)\} = \{1, 2\}$. So $c_1 \cup d_2$ is a
  colouring of $G$ that satisfies~\ref{c4}.
 \end{proof}
\Cref{thm:main2} immediately follows from \cref{l:c1,l:c2,l:c3,l:c4}.
 \end{proof}

 \section{Conclusion and open questions}
 \label{sec:c}

We collect here several  remarks and open questions.

\subsection{Fragile graphs have average degree less than 5}\label{subsec:av-deg-less-5}

As announced in the introduction, we recall the proof that every fragile graph $G$ on at least four vertices satisfies $|E(G)| \leq 2.5 |V(G)| -5$. 
When $G$ has 4 vertices, the inequality holds since the graph on~4 vertices and 6 edges is a complete graph and is 3-connected.  For the induction step, we  decompose $G$ into $G_1$ and $G_2$ as in the previous section. If $|V(G_1)| \leq 3$, then $G$ contains a vertex $x$ of degree at most~2.  Hence, $$|E(G)| \leq  |E(G\setminus x)| + 2 \leq 2.5|V(G\setminus x)| - 5 + 2 = 2.5(|V(G)| -1) -3 \leq 2.5|V(G)| -5.$$  We may therefore assume that $|V(G_1)| \geq 4$ and symmetrically $|V(G_2)| \geq 4$. Hence the induction hypothesis can be applied to both $G_1$ and $G_2$ so that the result follows from these inequalities: 


\begin{align*}
|E(G)| &\leq |E(G_1)| + |E(G_2)| \\
        &\leq 2.5|V(G_1)| - 5 + 2.5|V(G_2)| - 5 \\
         &= 2.5(|V(G_1)| + |V(G_2)|) - 10 \\
          &\leq 2.5 (|V(G)| +2) -10 \\
          &= 2.5 |V(G)| - 5.
\end{align*}

\subsection{Girth conditions}

It is easy to prove by induction that every fragile graph of girth at~least~4 on at~least~3 vertices satisfies $|E(G)| \leq 2 |V(G)| - 4$ (the proof is as in~\cref{subsec:av-deg-less-5}). This implies that every fragile graph with girth at least 4 contains a vertex of degree at most~3, so is \mbox{4-colourable}.  
We tried to improve this bound, but we instead found a fragile graph with girth~4 and chromatic number~4, as we now present.

\begin{figure}[h!]
  \begin{center}
    \includegraphics[width=10cm]{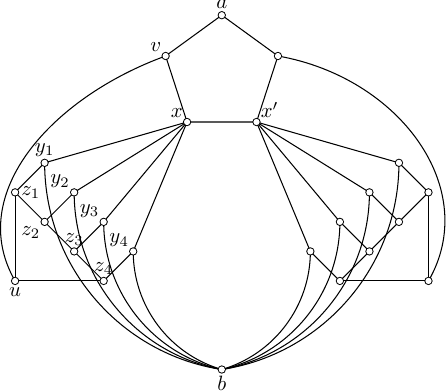}
  \end{center}
     \caption{The graph $G_1$.\label{f:G1}}
\end{figure}

Let $G_1$ be the graph represented in Figure~\ref{f:G1}.  It has girth~4 and is 2-degenerate; so in particular it is fragile and has chromatic number at~most~3. For all 3-colourings of~$G_1$, vertices $a$ and $b$ receive different colours.  Indeed, suppose for a contradiction that for some 3-colouring
of $G_1$, $a$ and $b$ receive the same colour, say colour~1.  Then, one of $x$ and $x'$, say $x$ up to symmetry, must receive a colour different from~1, say colour~2. So, the vertices $y_1$, \dots, $y_4$ must all receive the same colour, say colour~3. It follows that the vertices $z_1, \dots, z_4$ are coloured with colour~1 and~2 alternately. Hence, $u$ receives colour~3.  Now, $v$ has three neighbors, namely $a$, $x$ and $u$ that are coloured with colours~1, 2 and~3 respectively, a~contradiction.

\begin{figure}[h!]
  \begin{center}
    \includegraphics[width=14cm]{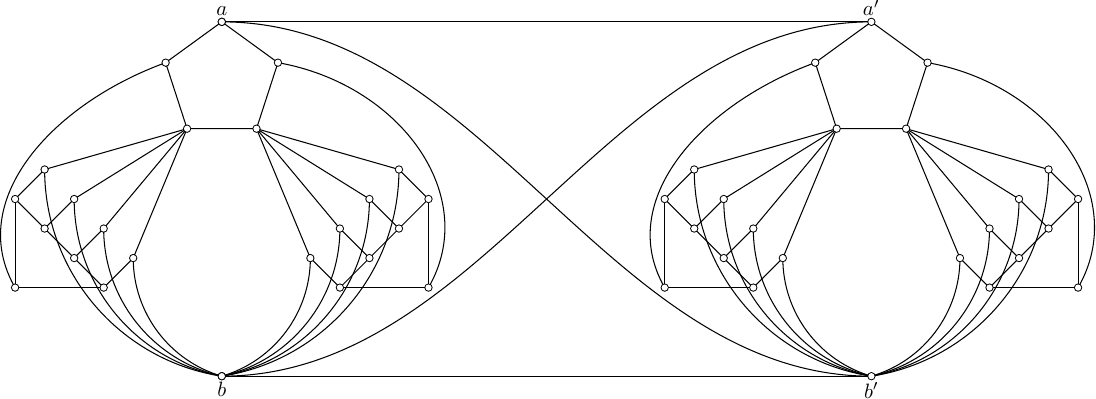}
  \end{center}
     \caption{The graph $G_2$.
\label{f:G2}}
\end{figure}

It follows that the triangle-free graph $G_2$ represented in Figure~\ref{f:G2} is not 3-colourable, but it is fragile since $\{a',b'\}$ is a~cutset, and $G_1$ is 2-degenerate even if two vertices adjacent to $a$ and $b$ are added.

\begin{figure}[h!]
  \begin{center}
    \includegraphics[height=7cm]{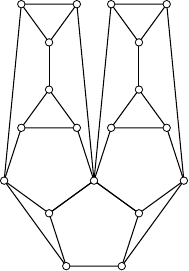}
  \end{center}
     \caption{A fragile graph with no cycle of length 4 and chromatic number~4. 
\label{f:noC4}}
\end{figure}

We could also obtain a fragile graph with no cycle of length~4 and chromatic number~4, see Figure~\ref{f:noC4}. 

This raises the following question: Is there a~finite girth that makes fragile graphs 3-colourable? A possible approach could be to prove that if the girth of a fragile graph is large enough, then the graph is 2-degenerate. But this approach fails because of the following construction. Consider an integer $g\geq 3$ and a connected cubic graph $G$ of girth $g$ (this exists, see for instance~\cite{ExJa:cage}).  Remove an edge $uv$ of $G$. This yields a 2-degenerate, and therefore fragile graph. Consider a copy $G'$ of $G\setminus uv$, with the vertices $u'$ and $v'$ corresponding to $u$ and $v$ respectively.  Now add an edge $uu'$ and an edge $vv'$. The obtained graph is fragile, cubic and has girth $g$.      


Trivially, a graph $G$ is fragile if and only if every subgraph $H$ of
$G$ is either on at most~3 vertices or admits a cutset of size at
most 2. In fragile graphs of girth at~least~4, one can further impose the cutset to be an independent set. 

\begin{lemma}
  A graph $G$ with girth at least~4 is fragile if and only if every subgraph $H$
  of~$G$ is either on at most~2 vertices or admits an independent cutset of size
  at most 2.   
\end{lemma}

\begin{proof}
  We prove the statement by induction on $|V(G)|$.
  The equivalence can be checked to hold on graphs of up to 3~vertices.  
If $|V(G)|\geq 4$, then since $G$ is not
  3-connected, it admits a cutset $S$ of size at most~2.  Suppose that $S$ is not
  independent, so $S=\{u, v\}$ and $uv\in E(G)$.  Let $C$ be a
  connected component of $G\setminus S$.  Since $G$ has girth at
  least~4, no vertex of $C$ is adjacent to both $u$ and $v$.  Hence,
  if $|C| = 1$, $G$ admits a cutset of size~1 (and therefore
  independent).  So we may assume that $|C| \geq 2$.  So, by the
  induction hypothesis, $G[S \cup C]$ admits an independent cutset
  $S'$.  It is easy to check that $S'$ is also a cutset of $G$.    
\end{proof}

\subsection{Algorithms}

By subdividing twice every edge of any graph $G$, a fragile graph $G'$ is obtained.  Poljak~\cite{poljak:74} proved that $\alpha(G') = \alpha(G) + |E(G)|$.  It follows that a polynomial-time algorithm that computes a maximum independent set for any fragile graph would yield a similar algorithm for all graphs.  This proves that computing a maximum independent set in a fragile graph is NP-hard. 

We also observe that, in $G'$, every edge $uv$ becomes a path $ux_{uv}y_{uv}v$.  Consider the graph $G''$ obtained from $G'$ by adding, for every vertex $x_{uv}$, a new vertex $x'_{uv}$ adjacent to $u$, $x_{uv}$ and $y_{uv}$.  It is easy to check that $G''$ is fragile and for all 3-colourings of $G''$ and all edges $uv$ of $G$, $u$ and $v$ have different colours (in $G''$).  It follows that if $G''$ is 3-colourable, then so is $G$.  Conversely it is easy to check that if $G$ is 3-colourable, so is $G''$. This proves that deciding whether a fragile graph is 3-colourable is NP-complete.  By the same kind of argument, we can prove that deciding whether a graph is 3-colourable stays NP-complete even when we restrict ourselves to fragile triangle-free graphs.  To see this, pick any graph $G$, remove all edges $uv$, and replace them by a copy of the graph $G_1$ from Figure~\ref{f:G1} with $a$ identified to $u$ and $b$ identified to $v$.  This yields a triangle-free fragile graph that is 3-colourable if and only if $G$ is 3-colourable. 

Our proof that every fragile graph is 4-colourable yields an
algorithm that actually computes a 4-colouring. A crude implementation of this algorithm would run in exponential time, but it is easy to turn it into a polynomial time algorithm by maintaining for each 2-tuples and 3-tuples $X$ of vertices of the input graph, a colouring satisfying the constraints \ref{c1}--\ref{c4} when applicable to $X$.

\section*{Acknowledgements}

\'Edouard Bonnet, Carl Feghali, St\'ephan Thomass\'e and Nicolas Trotignon  are partially supported by the French National Research Agency under research grant ANR DIGRAPHS ANR-19-CE48-0013-01, ANR Twin-width ANR-21-CE48-0014-01 and the LABEX MILYON  (ANR-10-LABX-0070) of Université de Lyon, within the program Investissements d’Avenir (ANR-11-IDEX-0007) operated by the French National Research Agency (ANR). 
Tung Nguyen and Paul Seymour are partially supported by AFOSR grants A9550-19-1-0187 and FA9550-22-1-0234, and by NSF grants DMS-1800053 and DMS-2154169.
Alex Scott is supported by EPSRC grant EP/X013642/1.
Part of this work was done when Nicolas Trotignon visited Princeton University with generous support of the H2020-MSCA-RISE project CoSP- GA No. 823748.


\end{document}